\documentclass[a4paper]{article}
\usepackage{amsmath}
\usepackage{amsthm}
\usepackage{amssymb}
\usepackage{graphicx}
\usepackage{authblk}
\usepackage{hyperref}
\usepackage{multirow}
\usepackage{cite}
\usepackage{comment}

\usepackage{caption}
\usepackage{subcaption}

\newtheorem{theorem}{Theorem}[section]

\newtheorem{conjecture}{Conjecture}

\newtheorem{remark}{Remark}
\newtheorem{example}{Example}
\newenvironment{vf}{\left\{\begin{array}{rcl}}{\end{array}\right.}
\DeclareMathOperator{\divergenceOperator}{div}
\DeclareMathOperator{\cycl}{Cycl}

\author[1]{Renato Huzak}
\author[2]{Goran Radunovi\'{c}}
\author[3]{Vesna \v{Z}upanovi\'{c}}

\affil[1]{Hasselt University, Campus Diepenbeek, Agoralaan Gebouw D, 3590 Diepenbeek, Belgium}
\affil[2]{University of Zagreb, Faculty of Science, Horvatovac 102a, 10000 Zagreb, Croatia}
\affil[3]{University of Zagreb, Faculty of Electrical Engineering and Computing, Unska 3, 10000 Zagreb, Croatia}

\title{Fractal analysis of slow-fast and regular systems: A survey of recent results and future perspectives}

\date{}

\begin{document}
\maketitle

\begin{abstract}

We survey recent developements in fractal analysis of regular and slow-fast dynamical systems using Minkowski dimension. Our focus is on spiral trajectories near monodromic limit periodic sets in regular systems and entry-exit sequences in slow-fast systems with degenerate singularities.

For regular systems, we recall connections between Minkowski dimension and cyclicity. Key results include fractal classifications of weak foci, degenerate foci, and polycycles, where dimensional relationships predict limit cycle birth.

For slow-fast systems, we survey the coordinate-free fractal methodology analyzing slow-fast Hopf points and canard cycles through slow divergence integrals and entry-exit sequences. The Minkowski dimension takes discrete values yielding upper bounds for the number of limit cycles without normal form transformations.

The fractal approach provides computational advantages, works with original coordinates, and reveals geometric structures underlying bifurcation phenomena. Applications span neuroscience, chemistry, population dynamics, and climate modeling. We also discuss extensions to piecewise smooth and three-dimensional systems.

\end{abstract}

\tableofcontents

\section{Introduction}

This paper contains an overview of recent developments in the fractal analysis of limit cycles and bifurcations in two-dimensional continuous dynamical systems. The study of fractal properties of dynamical systems goes back to Grassberger, Kaplan, Procaccia, Takens, Yorke, etc. 
Takens \cite{TakensFractal} showed that important information about deterministic dynamical systems (dimensions of attractors and entropy) can be obtained from time series produced by the systems without knowing the equations of motion. Using the notion of correlation dimension, Grassberger and Procaccia \cite{GrasProc} studied how to distinguish deterministic chaos and random noise in a dynamical system with a strange attractor. Kaplan and Yorke \cite{Kaplan} studied the dimension of invariant attractors using Lyapunov numbers. \v{Z}ubrini\'{c} and \v{Z}upanovi\'{c} \cite{Encikl} have provided a detailed overview of the
history of fractal dimensions and their applications in dynamical systems. See also \cite{paltak} and references therein. 

In this paper, we focus on a fractal analysis of limit periodic sets in families of planar vector fields that can generate limit cycles (isolated periodic orbits) after perturbation. For the analysis, we use the notions of Minkowski dimension (which
is equal to the box dimension \cite{Falconer}) and the corresponding Minkowski content (see Section \ref{subsection-Mink}). Limit cycles and oscillations appear in various problems in physics, biology, ecology, chemistry, medicine, meteorology, economics, etc. From a theoretical point of view, the determination of limit cycles in planar systems also plays a very important role (see Hilbert's 16th problem \cite{smale}). The cyclicity of a limit periodic set is the maximum number of limit cycles produced by that limit periodic set after perturbation (for a precise definition, we refer to e.g. \cite{RBook}). 

It is clear that the results on two-dimensional systems can be used in higher-dimensional systems by using invariant manifold reduction, etc.  

We distinguish two cases. In the first one (Sections \ref{subsection-regular-intro} and \ref{section-regular}) we deal with monodromic limit periodic sets, accumulated by spiral trajectories (foci, limit cycles, hyperbolic saddle-loops, hyperbolic 2-saddle polycycles). The Minkowski dimension of an arbitrary but fixed spiral trajectory or orbit of the associated Poincar\'{e} map is closely related to the number of limit cycles produced by such monodromic limit periodic sets. In the second case (Sections \ref{subsection-singular-intro} and \ref{section-slow-fast-case}) the limit periodic sets are degenerate (that is, they can be a slow-fast Hopf point or contain a curve of singular points). In this case, the Minkowski dimension of so-called entry-exit sequences in $\mathbb R$ is computed and related to cyclicity.

This paper provides an interesting survey of the connections between fractal geometry, qualitative theory of dynamical systems, including both regular and slow-fast systems, and numerical applications.

Note that the notion of Minkowski dimension is a natural one for studying fractal properties of orbits of dynamical systems. 
Namely, the classical Hausdorff dimension would not provide interesting results in our setting because of its countable stability, all of our objects would have trivial Hausdorff dimension and no information from it can be retrieved.
On the other hand, the failure of the Minkowski dimension to be countably stable is exactly what it makes it to be a useful tool for studying orbits of dynamical systems.


\subsection{The regular case}\label{subsection-regular-intro}
\smallskip

In the previous results of the authors, the Minkowski dimension of an orbit itself is considered  as well  as its dependence on a bifurcation parameter. 
Some well-known bifurcations are described in the new way using fractal analysis.
The Hopf and generalized Hopf bifurcations are studied in \cite{zuzu,ZupZub} (see also Section \ref{subsection-zuzu1}), while the saddle-node 
and period-doubling bifurcations of some discrete one-dimensional systems are studied in \cite{EZZ}. 
Under generic conditions the Minkowski dimension of an orbit at the bifurcation parameter are respectively 
$\frac43$, $\frac12$,  $\frac23$, for a spiral trajectory in the Hopf  bifurcation,  for one-dimensional orbits which undergo the saddle-node, and for the period doubling bifurcations.

Orbits tending to  hyperbolic singular or fixed points have trivial Minkowski dimension, that is $1$ for a  curve and $0$ for a sequence.
We can say that the Minkowski dimension reveals the quantity and the quality of objects born in a bifurcation, which is caused by the fact that the Minkowski dimension depends on the multiplicity of a singular or fixed point, see \cite{MRZ}.

Thus, the already  mentioned, but more degenerate bifurcations have orbits with larger Minkowski dimension. 
The results are based on the generalizations of formulas from \cite{tricot}.
As a simple example, we may consider a spiral $S$ defined in polar coordinates $r=\varphi^{-\alpha}$, $\varphi\ge\varphi_1>0$, where $\alpha\in(0,1]$, then according to the Tricot formula we obtain the Minkowski dimension $\dim_BS=  2/(1+\alpha)$. 

The Minkowski dimension is a useful tool for studying discrete and continuous dynamical systems, and the relation between them using the Poincar\' e map or the time-one map, see also \cite{BoxVesna,BoxDarko}. It is well known that $3$-dimensional systems exhibit interesting, chaotic behavior related to strange attractors which are fractal structured objects. For the Minkowski dimension of a class of $3$-dimensional systems studied in \cite{ZuZuR^3}, the dependence to coefficients of the system was found, contrary to one- and two-dimensional systems (see Section \ref{subsection-zuzu1}).

An interesting direction of the fractal analysis is related to oscillatory integrals, where the Minkowski dimension of a curve defined by the integral reveals the multiplicity of the phase function, see \cite{oscil}. The best known such type of integrals are Fresnel integrals and the related curve called the clothoid.

Our main goal is to give an overview of recent developments in the fractal analysis of degenerate foci (resp. polycycles) in planar analytic vector fields. We refer to Section \ref{subsection-degenerate} (resp. Section \ref{subsection-policycles-novo}).

 \medskip

\subsection{The slow-fast case}\label{subsection-singular-intro}
It is well-known that slow-fast systems play an important role in the study of dynamical systems, particularly in the analysis of periodic orbits and bifurcations. These systems are characterized by the presence of variables that evolve on significantly different time scales. Slow-fast systems arise in various fields, including neuroscience (Hodgkin-Huxley model \cite{Rubin2007}, bursting behaviors in neurons \cite{Saggio2017}, $\dots$), (bio)chemistry, e.g.,\cite{autocatalator,oxid,Segel1989}, etc.,
    population dynamics with predator-prey models, e.g., \cite{SuZhang,ChenLi,Yao2024}, etc.,
    climate modeling of glacial cycles \cite{Crucifix2012}, etc.
Slow-fast systems also play a crucial role in the theoretical study of limit cycles and bifurcations. 

Our main goal is to present some recent results \cite{BoxNovo,boxNovoNovo,BoxRenato,BoxDomagoj,BoxVlatko,BoxDarko} on fractal analysis of planar slow-fast systems which has proven to be a powerful tool for detecting the codimension of slow-fast (or singular) Hopf bifurcations in a coordinate-free way, finding the maximum number of limit cycles produced by canard cycles, etc. For the fractal analysis we use the notion of Minkowski dimension. 

The analysis of planar slow-fast systems frequently involves
studying small-amplitude limit cycles in the vicinity of a slow-fast
Hopf point (or a singular Hopf point). This is often called the birth of canards in the literature. We refer to Section \ref{section-slow-fast-case} for more details about slow-fast Hopf points. The codimension 1 slow-fast Hopf case
has been studied in \cite{KS} (a generalization of the Van der Pol system \cite{1996}). Slow-fast Hopf points of higher codimension in Li\'{e}nard vector fields have been treated in \cite{DRbirth}. Two important results in \cite{DRbirth} are: (a) finite upper bounds for the number of limit cycles in analytic
families and (b) finite upper bounds in smooth families with a finite codimension.

 In order to be able to use results from \cite{DRbirth,KS}, we first need to bring the slow-fast system into its normal form (see the classical Li\'{e}nard equation (11) in \cite{DRbirth} or the canonical form (3.3)  in \cite{KS}). This is sometimes very challenging from a computational point of view. In Section \ref{section-slow-fast-case}, we present a fractal method \cite{BoxNovo,BoxVlatko} for a direct
determination of upper bounds for the number of limit cycles near analytic slow-fast Hopf points (see Theorem \ref{theorem-Main1} and Remark \ref{remark-important} in Section \ref{subsection-Hopf}). By ``direct determination'' we mean working with the original coordinates (there is no need to use normal forms). The method consists in computing the Minkowski dimension of so-called entry-exit sequences in $\mathbb R$ associated with the slow-fast Hopf point (see Section \ref{subsection-SDI}). Then we obtain the upper bounds directly from the Minkowski dimension. There are simple formulas for computing the Minkowski dimension of the entry-exit sequences (for more details, we refer to \cite{BoxVlatko,BoxVesna} and Section \ref{subsection-Hopf}).

Section \ref{subsection-canard} is devoted to a fractal classification of canard cycles  that contain both attracting and repelling branches of a curve of singularities and a slow-fast Hopf point between them (see Figure \ref{fig-Canard}). Theorem \ref{theorem-Main2} connects the Minkowski dimension of appropriate entry-exit sequences to the number of limit cycles produced by such canard cycles, in a similar fashion to that used in Theorem \ref{theorem-Main1}.

A fractal analysis of canard cycles with two Hopf points can be found in \cite{BoxDomagoj}. We point out that the papers \cite{BoxRenato,BoxDomagoj,BoxVlatko} deal with bounded canard cycles. In \cite{boxNovoNovo}, one can find a fractal classification of slow-fast polynomial Li\'{e}nard equations near infinity (canard cycles are unbounded). 

For better readability, we present our results using a traditional standard form \eqref{eq-def-slowfast-1} in Section \ref{section-slow-fast-case}. We point out that the same results are true for planar slow-fast
systems given in non-standard form where a global separation of slow and fast variables is not possible (see \cite{BoxNovo,BoxVlatko}).

\subsection{Minkowski dimension and its properties} \label{subsection-Mink}

In this section, we briefly recall the definition of the Minkowski dimension \cite{Falconer,Goran}. Consider a bounded set $U$ in $\mathbb{R}^N$ and define its closed $\delta$-neighborhood as

\[
U_\delta = \{x \in \mathbb{R}^N \mid d(x,U) \leq \delta\}
\]
where $d$ is the Euclidean distance and by $|U_\delta|$ we denote the Lebesgue measure of $U_\delta$.

For any $s \geq 0$, the upper $s$-dimensional Minkowski content of $U$ is given as
\[
\mathcal{M}^{*s}(U) = \limsup_{\delta \to 0} \frac{|U_\delta|}{\delta^{N-s}}
\]
and similarly, the lower $s$-dimensional Minkowski content $\mathcal{M}_*^{s}(U)$ by replacing $\limsup$ with $\liminf$. When $\mathcal{M}_*^{s}(U)=\mathcal{M}^{*s}(U)$, the common value is called the $s$-dimensional Minkowski content of $U$, and denoted by $\mathcal{M}^{s}(U)$.

The lower and upper Minkowski (or box) dimensions of $U$ are defined as:

\[
\underline{\dim}_B U = \inf\{s \geq 0 \mid \mathcal{M}_*^s(U) = 0\}
\]
\[
\overline{\dim}_B U = \inf\{s \geq 0 \mid \mathcal{M}^{*s}(U) = 0\}
\]
When these dimensions coincide, we refer to it as the Minkowski dimension of $U$, denoted by $\dim_B U$.

A set $U$ is considered Minkowski nondegenerate if there exists a $d$ such that:
\[
0 < \mathcal{M}_*^d(U) \leq \mathcal{M}^{*d}(U) < \infty,
\]
in which case, $d = \dim_B U$, necessarily. We say that $U$ is Minkowski measurable if there exists $M^d(U)$ for some $d$ such that $M^d(U)\in (0,\infty)$.

A mapping $\Phi: U \subset \mathbb{R}^N \to \mathbb{R}^{N_1}$ is bi-Lipschitz if there exists a constant $\rho > 0$ such that:
\[
\rho\|x-y\| \leq \|\Phi(x)-\Phi(y)\| \leq \frac{1}{\rho} \|x-y\|
\]
for all $x,y \in U$.
The Minkowski dimensions are invariant under such maps:

\[
\underline{\dim}_B U = \underline{\dim}_B \Phi(U), \quad \overline{\dim}_B U = \overline{\dim}_B \Phi(U)
\]
Moreover, if $U$ is Minkowski nondegenerate, then $\Phi(U)$ is also Minkowski nondegenerate (see \cite{ZuZuR^3}). 

We also refer to \cite{Resman13} for some other properties of the Minkowski content.

\section{Fractal analysis of regular systems}\label{section-regular}

In Section \ref{subsection-zuzu1}, we recall some classical results in fractal analysis of weak foci and limit cycles. Section \ref{subsection-degenerate} (resp. Section \ref{subsection-policycles-novo}) is devoted to an overview of some recent results on the fractal properties of spiral trajectories near degenerate foci (resp. polycycles). 

We say that a spiral $r=f(\varphi)$ of focus type is {comparable} with the $\alpha$-{power spiral} $r=\varphi^{-\alpha}$ if $f(\varphi)/|\varphi|^{-\alpha}\in[ A, B]$ for some positive constants $ A$ and $ B$ and for all $\varphi\in [1,\infty)$ (resp. $\varphi\in (-\infty,-1]$) if the spiral has positive (resp. negative) orientation. Similarly, a spiral $r=f(\varphi)$ of focus type is said to be comparable with the exponential spiral $r=e^{-\beta\varphi}$ if $f(\varphi)/e^{-\beta\varphi}\in[ A, B]$, where $ A$ and $ B$ are positive constants, $\beta$ is a positive (resp. negative) constant and $\varphi\in [0,\infty)$ (resp. $\varphi\in (-\infty,0]$) if the spiral has positive (resp. negative) orientation.

In a similar fashion, we define comparability with $r=a\pm\varphi^{-\alpha}$ and with $r=a\pm e^{-\beta\varphi}$ where $r=a$ corresponds to a limit cycle ($a>0$). For more details, we refer the reader to \cite{zuzu}.

In Theorems \ref{zuzufocus}--\ref{ThmFocus3} (resp. Theorem \ref{ThmFocus4}) below, $(r, \varphi)$  are the standard (resp. generalized) polar coordinates. 

\subsection{Fractal analysis of  weak foci}\label{subsection-zuzu1}

We consider a standard generic Hopf-Takens bifurcation \cite{takens} in polar coordinates
\begin{equation}\label{polar}
                \begin{cases}
                    \dot{r} = r\left(r^{2l} + \sum_{i=0}^{l-1}a_ir^{2i}\right)\\
                    \dot{\varphi} = 1,
                \end{cases}
            \end{equation}
with parameters $a_0, \dots, a_{l-1}\in \mathbb R$. 

\begin{theorem}[Fractality of Hopf-Takens bifurcation  - focus case \cite{zuzu}]\label{zuzufocus}
                Let $S$ be part of a trajectory of the system (\ref{polar}) near the origin. 
                \begin{enumerate}
                    \item If $a_0\neq 0$, then the spiral $S$ is comparable
with the exponential spiral $r = e^{a_0\varphi}$ and ${\dim}_B\, S = 1$.
                    \item If $a_0 = a_1 = \ldots = a_{k-1} = 0, a_k \neq 0$ for $1\le k\le l$ ($a_l=1)$, then $S$ is comparable with the spiral $r = \varphi^{-1/2k}$ of power type and ${\dim}_B\, S = \frac{4k}{2k+1}$.
                \end{enumerate}
                In both cases $S$ is Minkowski measurable.
            \end{theorem}

For $a_0 \neq 0$, the focus is hyperbolic, so the trajectory tending to the origin is exponential and has finite length. On the other hand, in the case $a_0 = 0$, the focus is weak with pure imaginary eigenvalues. A trajectory near the origin is a non-rectifiable curve, and the Minkowski dimension is obtained using Tricot’s formula from \cite{tricot} for spirals of the form $r = \varphi^{-\alpha}$, $0 < \alpha \le 1$, $\varphi>0$,  along with comparison criteria obtained in \cite{zuzu}.
\smallskip

\begin{figure}[htb]
	\begin{center}
		\includegraphics[width=9.3cm]{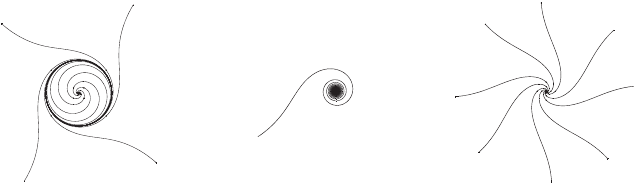}
		{\footnotesize 
        \put(-236,-10){$a_0<0$}
         \put(-136,-10){$a_0=0$}
          \put(-39,-10){$a_0>0$}
}
         \end{center}
	\caption{The standard Hopf bifurcation $\dot r=r(r^2+a_0),\dot\varphi=1$, for $l=1$ in \eqref{polar}. For $a_0<0$, spiral trajectories near the strong focus $r=0$ or near the limit cycle $r=\sqrt{-a_0}$ are of exponential type and hence the Minkowski dimension is $1$ (see also Theorem \ref{zuzulc}). When $a_0=0$, the Minkowski dimension of spiral trajectories near the weak focus $r=0$ is $4/3$. For $a_0>0$, all trajectories near the strong focus $r=0$ have Minkowski dimension equal to $1$.}
	\label{fig-HopfTakens}
\end{figure}

   Notice that system (\ref{polar}) is a normal form for a weak focus, see \cite{takens}. Furthermore, for codimension $l=1$ the system undergoes the Hopf bifurcation while for $l\ge 2$ the degenerate Hopf bifurcation or the Hopf-Takens bifurcation occurs. The codimension $l$ is the number of limit cycles that can be produced in the system (\ref{polar}). By using Theorem \ref{zuzufocus} a correspondence is established between the Minkowski dimension and  the cyclicity of the origin $r=0$ in (\ref{polar}).
   Standard Hopf bifurcation appears for $l=1$,  the Minkowski dimension of a spiral trajectory is equal to $4/3$, and $1$ limit cycle is produced in the bifurcation (see the example in Figure \ref{fig-HopfTakens} taken from \cite{zuzu}).

\smallskip

A standard method for studying number of limit cycles includes the Poincar\' e map, also called the first return map. The number of limit cycles corresponds to the number of fixed points of the Poincar\' e map near focus, limit  cycle or polycycle. 
The result from Theorem \ref{zuzufocus} could be obtained from the Minkowski dimension of a discrete orbit generated by the Poincar\' e map, see \cite{ZupZub} for such approach.

The following theorem deals with fractal analysis of limit cycles in (\ref{polar}).

\begin{theorem}[Fractality of Hopf-Takens bifurcation  - limit cycle case \cite{zuzu}]\label{zuzulc}

            Let the system (\ref{polar})            
        have a limit cycle $r = a$ of multiplicity $m$, $1\leq m\leq l$. Let $S_1$ and $S_2$ be spiral trajectories of this system near the limit cycle from the outside and inside respectively. 
        Then ${\dim}_B\, S_1 = {\dim}_B\, S_2 = 2 - \frac{1}{m}$.  The trajectories $S_1$ and $S_2$ are
comparable with exponential spirals $r = a \pm e^{-\beta\varphi}$, $\beta\ne 0$, of limit cycle type when $m = 1$, and with power spirals $r = a \pm \varphi^{-1/(m-1)}$ when $m > 1$.
        \end{theorem}

\smallskip
We refer to \cite[pg. 55]{RBook} for a definition of the notion of multiplicity of limit cycles.
The idea of the proof is analogous as for Theorem \ref{zuzufocus}. Furthermore, from a  limit cycle of  multiplicity $m=2$ another limit cycle could be produced. The Poincar\' e map near the limit cycle at a transversal,   undergoes generic saddle-node bifurcation. According to Theorem \ref{zuzulc} a spiral tending to the limit cycle has the Minkowski dimension $3/2$, while the dimension of the corresponding discrete orbit is equal to $1/2$, see \cite{ZupZub}. 
Theorems \ref{zuzufocus} and \ref{zuzulc} have also been generalized to the case of Hopf-Takens bifurcation at infinity \cite{GoranInf}.

\smallskip

\begin{figure}
\centering
\begin{subfigure}{.5\textwidth}
  \centering
  \includegraphics[width=.6\linewidth,height=\linewidth,angle=270,trim=280pt 280pt 280pt 280pt,clip]{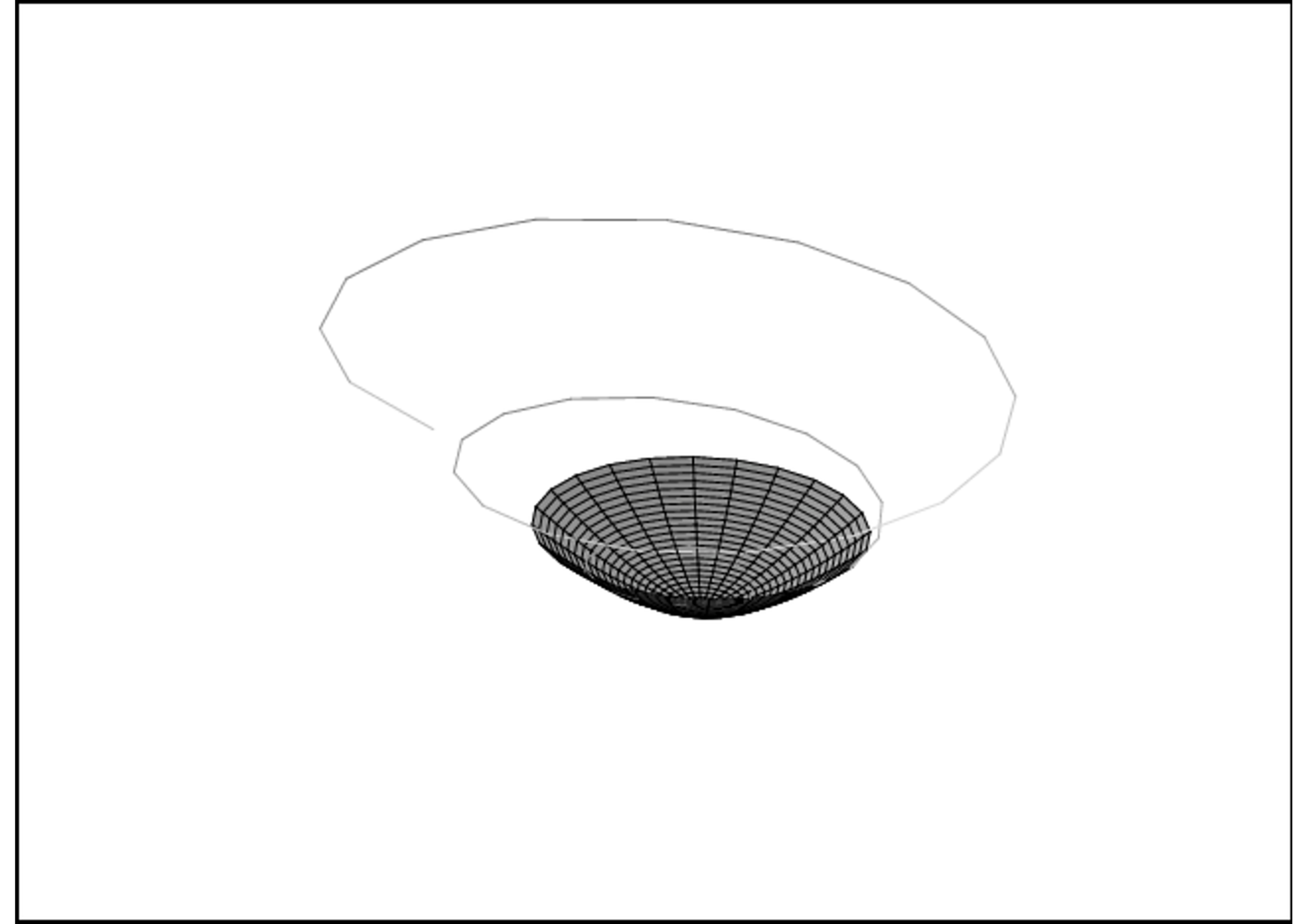}
  \caption{Lipschitzian spiral}
  \label{fig:sub1}
\end{subfigure}%
\begin{subfigure}{.5\textwidth}
  \centering
  \includegraphics[width=0.6\linewidth, height=\linewidth,angle=270,trim=250pt 280pt 280pt 280pt,clip]{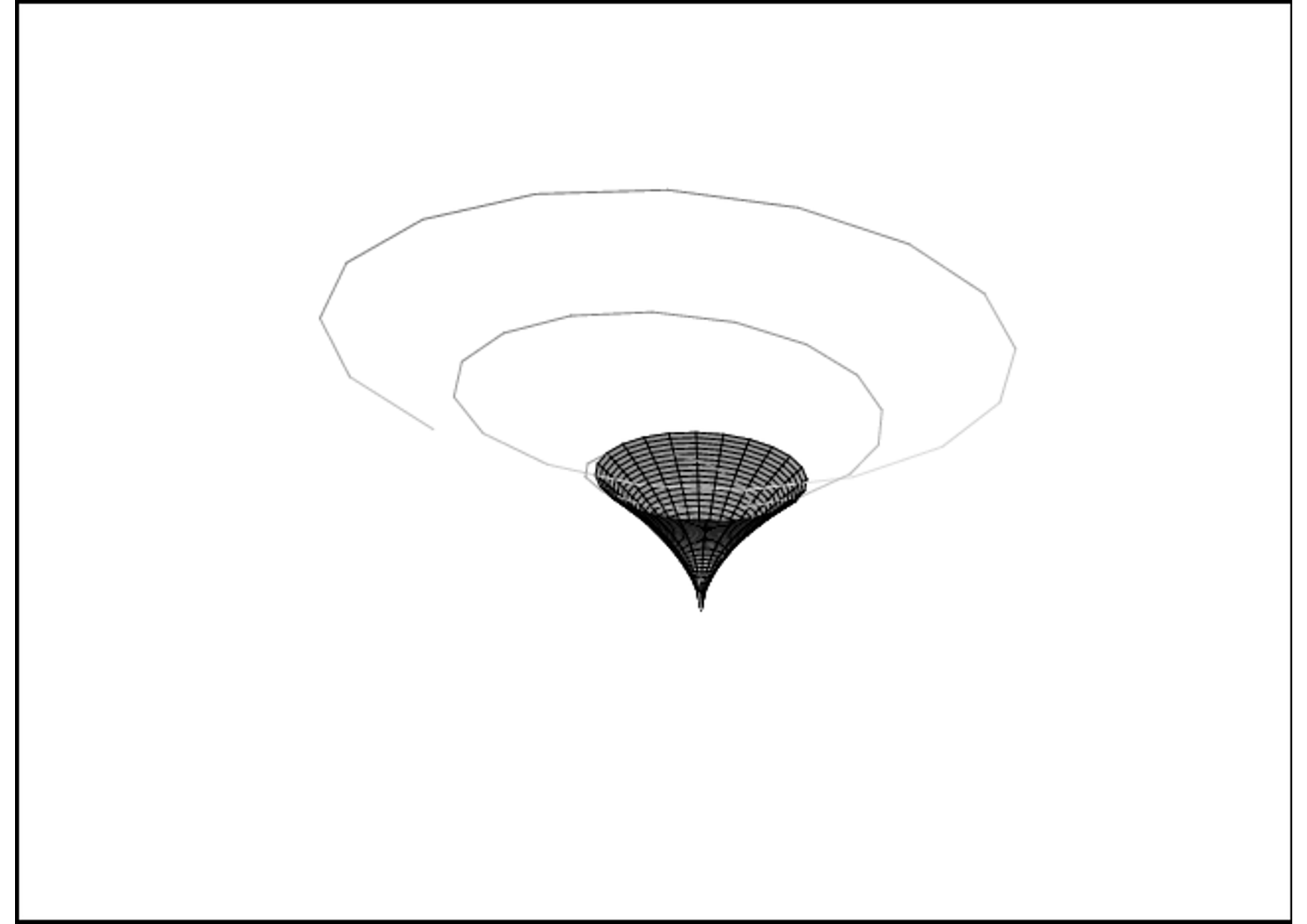}
  \caption{H\"olderian spiral}
  \label{fig:sub2}
\end{subfigure}
\caption{Spiral trajectories in $\mathbb R^3$.}
\label{fig:test}
\end{figure}

The results could be applied to $3$-dimensional systems having a linear part 
with pure imaginary pair and a simple zero eigenvalues,
see \cite{ZuZuR^3}.
Trajectories of the $3$-dimensional systems are considered in a two-dimensional surface. Depending on whether the surface is of Lipschitz or of H\"older type near the limit set, we distinguish 
the following two cases: Lipschitzian spirals  and H\"olderian spirals (see Fig. \ref{fig:test}). 
Since the Minkowski dimension is not affected by bi-Lipschitz map it is proved that the Minkowski dimension of a trajectory is equal to the dimension of the projection in the plane, which is weak focus studied in \cite{zuzu}.
The H\"olderian case is studied separately. 

In $\mathbb R^3 $ a new phenomenon showed up, that is, spiral trajectories
can be sensitive to coefficients of a system.

\begin{example} [Spiral trajectories in $\mathbb{R}^3$ \cite{ZuZuR^3}]

The solution of the system  
\begin{equation}
  \begin{cases}
\dot r=a_1rz \\
\dot \varphi =1 \\
\dot z=b_2z^2,
\end{cases}
\end{equation}
with parameters  $a_1, b_2\in \mathbb R$ is a family of spirals 

$$
r=C_1(-b_2t+C_3)^{-a_1/b_2},
\varphi=t+C_2,
z=\frac1{-b_2t+C_3}.
$$
Spirals $S$ near the origin are contained in the surface
$z=C\cdot r^{b_2/a_1}$.

(i) (Lipschitzian spirals) If $a_1/b_2\in(0,1]$ then $\dim_BS=\frac{2}{1+a_1/b_2}$.

(ii) (H\"olderian spirals) If $a_1/b_2>1$ then $\dim_BS=1$.

(iii) If $a_1/b_2<0$ then the origin is not an accumulation point for $S$.

\end{example}

\subsection{Fractal analysis of  degenerate foci }\label{subsection-degenerate}

   The weak focus studied in the previous section is the simplest focus case. 
      A class of singularities of focus type with two zero eigenvalues, nilpotent or more degenerate, is a generalization of the weak focus case. 
      The nilpotent focus has been studied from the fractal point of view in \cite{BoxVesna}, as well as the application of the fractal results to the cyclicity of a family of nilpotent foci by using results from \cite{hanrom}.

      It is well known that the blowing-up procedure leads to a decomposition of a singularity to elementary singularities (hyperbolic and semihyperbolic).
      The simplest case, the weak focus has no singularities obtained in the blowing-up procedure. The consequence is that the Poincar\' e map  at each transversal has the same asymptotics.
      For nilpotent and degenerate foci this is not the case. A nilpotent focus has one characteristic transversal curve, so the cyclicity could be read from the Minkowski dimension at the curve,  \cite{BoxVesna}.

      Furthermore, there are degenerate foci having no singularities  after the blowing-up procedure. 
   Degenerate spirals appear in complex swirling flows, see, e.g.\ \cite{fraser}.   Most naturally occurring spirals are developed in systems with inherent asymmetry.

   We consider the following systems \cite{BoxDarko} near the origin:

   \begin{equation}\label{system-deg-foci}
                \begin{cases}
                    \dot{x} = -ny^{2n-1}\pm nx^m y^{n-1} (x^{2m}+y^{2n})^k\\
                    \dot{y} = mx^{2m-1}\pm m x^{m-1} y^n (x^{2m}+y^{2n})^k,
                \end{cases}
            \end{equation}
            with $m,n\in\mathbb{N}\setminus\{0\}$, $m+n>2$, $m$ and $n$ odd and $k\in\mathbb N$. If the sign in \eqref{system-deg-foci} is negative (resp. positive), we have
a stable (resp. unstable) focus at the origin.  

\begin{theorem} [Fractality of the degenerate focus; $m=n$ case \cite{BoxDarko}]\label{ThmFocus3}
   Suppose that $m = n$ in \eqref{system-deg-foci}. Let ${S}$ be a spiral trajectory of \eqref{system-deg-foci} near the origin. Then the following statements are true.
\begin{enumerate}
    \item If $k=0$, then the spiral ${S}$ is comparable with the exponential spiral $r=e^{\pm\varphi/n}$,
    $\dim_B{S}=1$, and ${S}$ is Minkowski measurable.
    \item If $k>0$, then the spiral ${S}$ is comparable with the power spiral $r=\varphi^{-1/2nk}$,
$\dim_B{S}=2-\frac{2}{1+2kn}$,
and ${S}$ is Minkowski nondegenerate.
\end{enumerate}
\end{theorem}

In the rest of this section we focus on the case $m\ne n$.
We use the $(n,m)$--polar coordinates $$(x,y)=(r^n\text{Cs}(\varphi),r^m \text{Sn}(\varphi)),$$ where $\text{Cs}(\varphi)$ and $\text{Sn}(\varphi)$ are a generalization of $\cos\varphi$ and $\sin\varphi$, and
\[\frac{d}{d\varphi}{\text{Cs}}(\varphi)=-n \text{Sn}^{2n-1}(\varphi), \ \frac{d}{d\varphi}{\text{Sn}}(\varphi)=m \text{Cs}^{2m-1}(\varphi)\]
with $(\text{Cs}(0),\text{Sn}(0))=(1,0)$. It is not difficult to see that $\text{Cs}^{2m}(\varphi)+\text{Sn}^{2n}(\varphi)=1$, $\text{Cs}(\varphi)$ is even, $\text{Sn}(\varphi)$ is odd, and both are $T$-periodic, with 
\[T=\frac{2}{mn}\frac{\Gamma(\frac{1}{2m})\Gamma(\frac{1}{2n})}{\Gamma(\frac{1}{2m}+\frac{1}{2n})},\]
where $\Gamma$ is the gamma function (see \cite{BoxDarko}). System \eqref{system-deg-foci} becomes
\begin{equation}\label{spiralll}
\frac{dr}{d\varphi}=\pm \text{Sn}^{n-1}(\varphi)\text{Cs}^{m-1}(\varphi)r^{2mnk+1},\nonumber
\end{equation}
upon division of $\dot{r}$ by $\dot{\varphi}$. We have

\begin{theorem}[Fractality of the degenerate focus; $m\neq n$ case \cite{BoxDarko}] \label{ThmFocus4} Let $T$ be the period of the functions $\text{Cs}(\varphi)$ and $\text{Sn}(\varphi)$. Then the following statements are true for a spiral trajectory ${S}$ of \eqref{system-deg-foci} near the origin and in the case when $m\neq n$.
\begin{enumerate}
    \item If $k=0$, then the spiral ${S}$ is comparable with the exponential spiral $r=e^{\pm \frac{2\pi}{T m n}\varphi}$,
    $\dim_B{S}=1$, and ${S}$ is Minkowski measurable.
    \item If $k>0$, then the spiral ${S}$ is comparable with the power spiral $r=\varphi^{-1/2mnk}$.
\end{enumerate}
\end{theorem}

Based on numerical observations (see Table \ref{tablica_dims_2}), in \cite{BoxDarko} we proposed the following conjecture:
\begin{conjecture}
\textit{Let $k$ be a positive integer, $m> n$ and let $m,n$ be odd. Then any spiral trajectory of \eqref{system-deg-foci}
  near the origin $(x,y)=(0,0)$ has the Minkowski dimension $2-\frac{1+\frac{n}{m}}{1+2kn}$. }

\end{conjecture}

\begin{table}[htb]
\begin{center}
\begin{tabular}{|c|c|c|c|c|c|}
\hline
$m$ & $n$ & $k$ & \multicolumn{2}{c|}{\begin{tabular}[c]{@{}c@{}}conjectured\\ dimension\end{tabular}} & \begin{tabular}[c]{@{}c@{}}numerical\\ dimension\end{tabular} \\ \hline
5   & 3   & 2   & $122/65$                                   & $1.87692$                               & $1.87287$                                                     \\ \hline
11  & 3   & 2   & $272/143$                                  & $1.90210$                               & $1.89615$                                                     \\ \hline
21  & 3   & 2   & $174/91$                                   & $1.91209$                               & $1.90574$                                                     \\ \hline
21  & 11  & 2   & $1858/945$                                 & $1.96614$                               & $1.96561$                                                     \\ \hline
5   & 3   & 11  & $662/335$                                  & $1.97612$                               & $1.97581$                                                     \\ \hline
11  & 3   & 11  & $1460/737$                                 & $1.98100$                               & $1.98063$                                                     \\ \hline
21  & 3   & 11  & $930/469$                                  & $1.98294$                               & $1.98255$                                                     \\ \hline
21  & 11  & 11  & $10174/5103$                               & $1.99373$                               & $1.99355$                                                     \\ \hline
\end{tabular}
\caption{Conjectured and numerical box dimensions of $S$ in \eqref{system-deg-foci}. The numerical estimates are computed in Wolfram Mathematica, version 12, based on a decomposition of spiral trajectories $S$ using circular sectors and on the concept of nucleus and tail. For more details, we refer the reader to \cite{BoxDarko}. }
\label{tablica_dims_2}
\end{center}
\end{table}

\subsection{Fractal analysis of polycycles}\label{subsection-loop}\label{subsection-policycles-novo}

\begin{figure}[htb]
	\begin{center}
		\includegraphics[width=9.3cm,height=2.9cm]{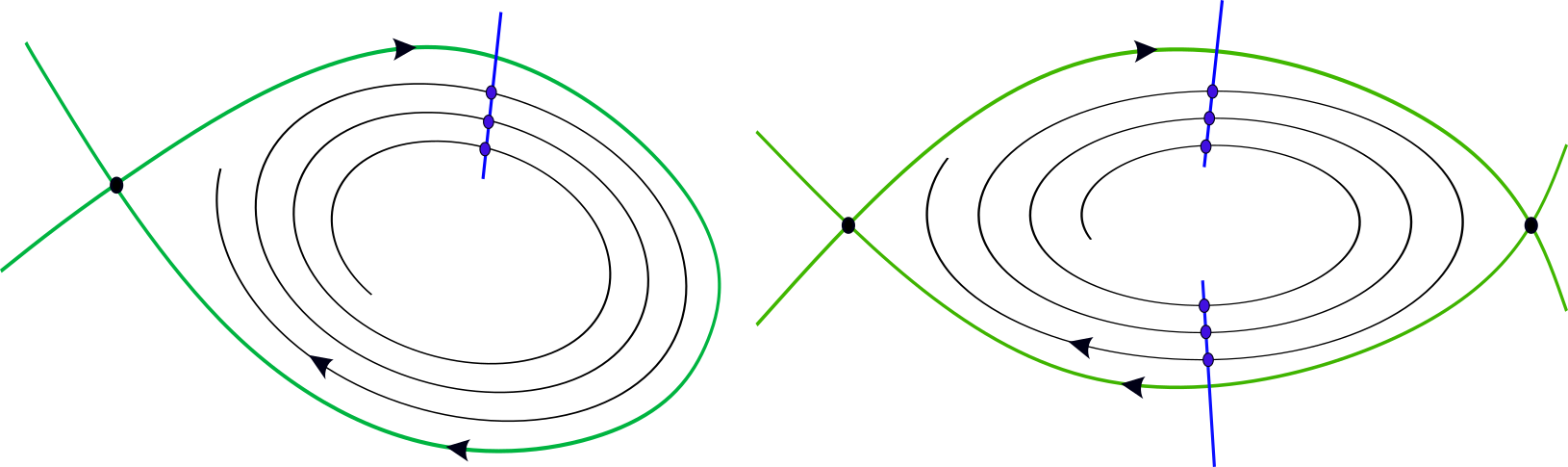}
		{\footnotesize 
        \put(-200,-10){$(a)$}
        \put(-186,83){$t_1$}
        \put(-206,40){$S$}
        \put(-70,-10){$(b)$}
        \put(-67,81){$t_1$}
        \put(-58,4){$t_2$}
        \put(-80,40){$S$}
}
         \end{center}
	\caption{(a) A hyperbolic saddle-loop ($\tilde N=1)$. (b) A hyperbolic 2-saddle polycycle ($\tilde N=2$).}
	\label{fig-polycyc}
\end{figure}
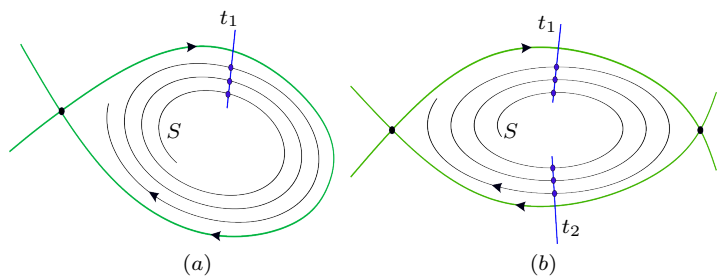

In this section we present some new results \cite{vlatko,crhure} on the fractal analysis of spiral trajectories near polycycles in planar analytic vector fields. Polycycles consist of a finite number of hyperbolic saddles and regular connections between them (see Figure \ref{fig-polycyc}). The Poincar\' e map to transversals to a
polycycle need not be analytic. 


We start with the following general result.
 \begin{theorem}[Fractality of a hyperbolic polycycle \cite{crhure}]\label{theorem-general-poli}
        Consider a polycycle of an analytic vector field, with $\tilde N$ hyperbolic saddles as vertices, $\tilde N=1,2,\dots$. 
        Let $S$ be a spiral trajectory that accumulates to the polycycle. Let $t_1, t_2, ...., t_{\tilde N}$ be transversals to $($all$)$ regular sides of the polycycle (see Figure \ref{fig-polycyc}). By $Y_i,\, i\in\{1,2,...,\tilde N\}$, we denote the intersections of $S$ with $t_i$.
        Then,
        \[ \dim_B S = 1 + \max \left\{ \dim_B Y_i \colon i\in\{1,2,...,\tilde N\} \right\}. \]
    \end{theorem}

If $p$ is a hyperbolic saddle with eigenvalues $\eta_-<0<\eta_+$, then we call $r = -\frac{\eta_-}{\eta_+}$ the hyperbolicity ratio of $p$. Theorem \ref{theorem-looop} below shows that the correspondence
between the notion of codimension \cite{RBook} of the saddle-loop and the Minkowski dimension of spiral trajectories is 2–1. Following \cite{RBook}, the codimension of the saddle-loop
corresponds to its cyclicity in generic unfoldings.

\begin{theorem}[Fractality of a saddle-loop \cite{crhure}]\label{theorem-looop}
                The Minkowski dimension of a spiral trajectory $S$ in an analytic vector field that has a finite-codimension saddle-loop as its $\alpha/\omega$-limit set depends only on the codimension of the saddle-loop.
                More precisely, if $k\geq 1$ is the codimension of the saddle-loop, then:
                \[ \dim_B S = \begin{cases}
                2 - \frac{2}{k},& k \text{ even,}\\
                2 - \frac{2}{k+1},& k \text{ odd}. \end{cases} \]
        \end{theorem}

The following result implies that we can read the known cyclicity bounds\cite{Mour3} for hyperbolic 2-saddle
polycycles using only Minkowski dimensions of spiral trajectories and associated sequences of intersection points. 
   
\begin{theorem} [Fractality of hyperbolic 2-saddle polycycles \cite{vlatko}]
                Consider a hyperbolic 2-saddle polycycle of an analytic vector field $X_0$ and denote by $r_1$ and $r_2$ ratios of hyperbolicity of polycycle saddles. Assume that $r_1r_2 \ne  1$ or $r_1 \notin \mathbb{Q}$. Let $S$ be its one accumulating spiral trajectory. Then the following statements are true. 
                \begin{enumerate}
                \item If $\dim_B S=1$, then the cyclicity of the polycycle in $C^\infty$ unfoldings of $X_0$ is at most $3$.
                \item If $d:=\dim_B S\in(1,2)$, then the upper bound on the cyclicity of the polycycle in $C^\infty$ unfoldings of $X_0$ is given by 
                
                \begin{equation}
            \left\lfloor 3 + (1+r)\frac{d-1}{2-d} \nonumber\right\rfloor,\end{equation} where
                $$
                r=\min\left\{\frac{d_2(1-d_1)}{d_1(1-d_2)},\frac{d_1(1-d_2)}{d_2(1-d_1)}\right\},
                $$
                and $d_1\in(0,1)$ and $d_2\in(0,1)$ are Minkowski dimensions of sequences obtained as intersections of spiral $S$ with transversals to the two heteroclinic connections. Note also that $d=1+\max\{d_1,d_2\}$ (see Theorem \ref{theorem-general-poli}).
                \end{enumerate}
                \end{theorem}

We refer the reader to \cite{vlatko,Mour3} for a definition of $C^\infty$ unfoldings.

\section{Fractal analysis of slow-fast systems}\label{section-slow-fast-case}

We consider a two-dimensional slow-fast system
\begin{equation}\label{eq-def-slowfast-1}
X_{\epsilon,\lambda}: \ \left\{
\begin{array}{rcl}
\dot{x} & = & f(x,y,\epsilon,\lambda), \\
\dot{y} & = & \epsilon g(x,y,\epsilon,\lambda),
\end{array}
\right.
\end{equation}
where $0<\epsilon\ll 1$ represents the singular perturbation parameter, $\lambda\sim\lambda_0\in\mathbb R^m$ and $f$ and $g$ are $C^\infty$-smooth functions. 

We assume that the fast subsystem $X_{0,\lambda}$ has a curve of singularities $$\mathcal{C}_{\lambda}:=\{(x,y)\in\mathbb R^2\mid f(x,y,0,\lambda)=0\},$$
for all $\lambda$ kept close to $\lambda_0$. A singularity $p\in \mathcal C_\lambda$ is normally hyperbolic if $\frac{\partial f}{\partial x}(p,0,\lambda)\neq 0$ (attracting when $\frac{\partial f}{\partial x}(p,0,\lambda)< 0$ or repelling when $\frac{\partial f}{\partial x}(p,0,\lambda)> 0$). If the normal hyperbolicity is lost at some point $p\in \mathcal C_\lambda$, i.e.,
$$\frac{\partial f}{\partial x}(p,0,\lambda)=0,$$
then we call the singularity $p$ a contact point (the curve of singularities $\mathcal{C}_{\lambda}$ has a contact with the horizontal fast orbits of $X_{0,\lambda}$ at $p$). We refer to Figure \ref{fig-Hopf}. The contact point $p$ is nilpotent if
$$\frac{\partial f}{\partial y}(p,0,\lambda)\neq0.$$

Finally, a nilpotent contact point $p\in \mathcal C_\lambda$ is called a slow-fast Hopf point \cite[Definition 2.4]{DDR-book-SF} if 
\begin{equation}\label{eq-def-hopf-turning-point}
\begin{split} 
  g(p,0,\lambda)=0, \ \ \frac{\partial^{2} f}{\partial x^{2}}(p,0,\lambda) \neq 0, \ \ \frac{\partial g}{\partial x}(p,0,\lambda)\cdot \frac{\partial f}{\partial y}(p,0,\lambda) < 0.\nonumber
\end{split}
\end{equation}
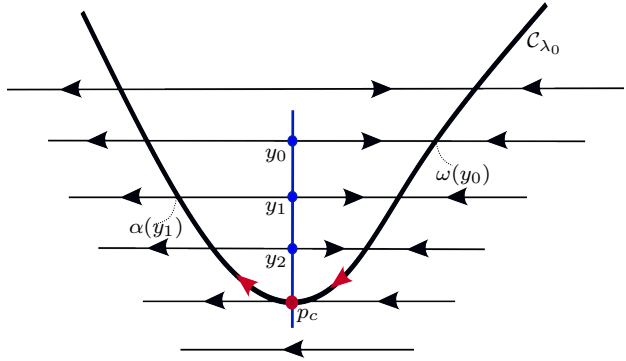
\begin{figure}[htb]
	\begin{center}
		\includegraphics[width=8.1cm,height=4.7cm]{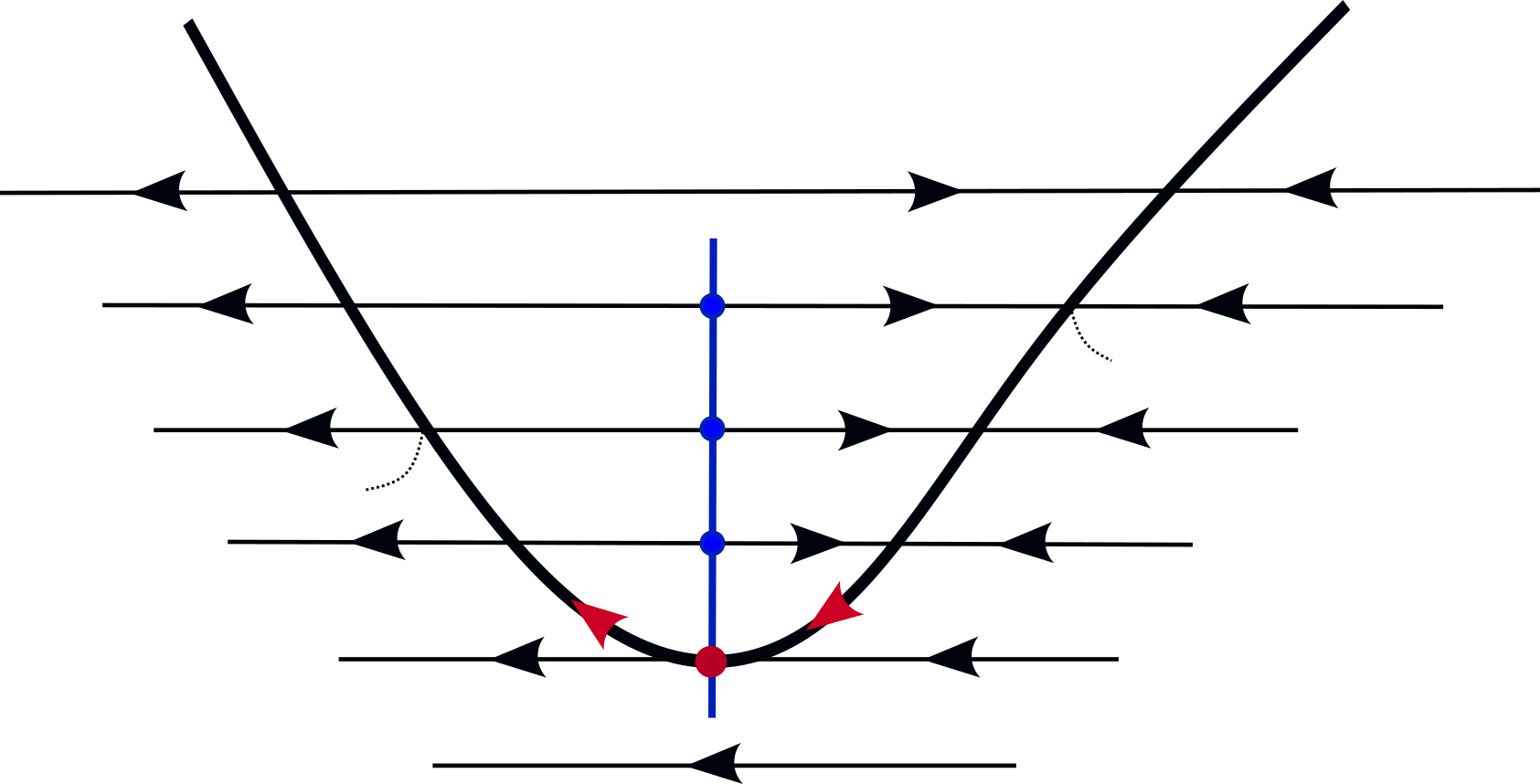}
		{\footnotesize \put(-122,15){$p_c$}
        \put(-36,117){$\mathcal C_{\lambda_0}$}
        \put(-70,67){$\omega(y_0)$}
        \put(-185,47){$\alpha(y_1)$}
        \put(-134,75){$y_0$}
        \put(-134,55){$y_1$}
        \put(-134,35){$y_2$}
}
         \end{center}
	\caption{Dynamics of $X_{0,\lambda_0}$ near the slow-fast Hopf point $p_c$ when  
   $\frac{\partial^{2} f}{\partial x^{2}}(p_c,0,\lambda_0)<0$, $\frac{\partial f}{\partial y}(p_c,0,\lambda_0)>0$ and $\frac{\partial g}{\partial x}(p_c,0,\lambda_0)<0$. $(y_k)_{k\ge 0}$ is an entry-exit sequence.   }
	\label{fig-Hopf}
\end{figure}
Now, we assume the following. \\
\\
\textbf{Assumption 1} \textit{(Slow-fast Hopf point).} The slow-fast family $X_{\epsilon,\lambda}$ defined by \eqref{eq-def-slowfast-1} has a slow-fast Hopf point $p_c=(x_c,y_c)$ at $\lambda_0$.\\
\\
From Assumption 1 and the Implicit Function Theorem it follows that $p_c$ is a contact point of Morse type (i.e. of
contact order 2) and separates the normally attracting branch and the normally repelling branch of $\mathcal C_{\lambda_0}$. Moreover, near $p_c$, the curve of singularities $\mathcal C_{\lambda_0}$ is concave up (resp. concave down) when $\frac{\partial f}{\partial y}(p_c,0,\lambda_0)$ and $\frac{\partial^{2} f}{\partial x^{2}}(p_c,0,\lambda_0)$ have opposite signs (resp. the same sign).

Since fractal sequences (often called entry-exit sequences) tending monotonically to $p_c$ play a very important role in this section, first we introduce the concept of a slow divergence integral of $X_{\epsilon,\lambda}$ and entry-exit function (Section \ref{subsection-SDI}). Then we give a complete fractal classification of $p_c$ (Section \ref{subsection-Hopf}). For the classification, we used the Minkowski dimension of the entry-exit sequences. When $X_{\epsilon,\lambda}$ is analytic, we connect the Minkowski dimension to upper bounds for the number of limit cycles of $X_{\epsilon,\lambda}$ produced by $p_c$.

In Section \ref{subsection-canard}, we classify the bounded canard cycles of $X_{\epsilon,\lambda}$. Again, we use the Minkowski dimension of appropriate entry-exit sequences and connect the dimension to the number of limit cycles of $X_{\epsilon,\lambda}$ (Hausdorff close to such canard cycles).

\subsection{Slow divergence integral and entry-exit sequences}\label{subsection-SDI}

We  write the slow-fast family \eqref{eq-def-slowfast-1} as
$$X_{\epsilon,\lambda}=f(x,y,0,\lambda)\frac{\partial}{\partial x}+\epsilon Q_\lambda+O(\epsilon^2),$$
where $f$ is introduced in \eqref{eq-def-slowfast-1} and $Q_\lambda$ is a smooth family of vector fields. 
We define the notion of slow vector field along normally hyperbolic segments of $\mathcal C_\lambda$. Our definition follows that of \cite[Chapter 3]{DDR-book-SF}. If $p=(x,y)\in \mathcal C_\lambda$ is a normally hyperbolic singularity, then $\bar Q_\lambda(p)\in T_{p}\mathcal{C}_{\lambda}$ denotes the projection of $Q_\lambda(p)$ on $T_{p}\mathcal{C}_{\lambda}$ in the direction of the fast horizontal orbits. $T_{p}\mathcal{C}_{\lambda}$ is the eigenspace of the zero eigenvalue of $X_{0,\lambda}$'s linear part at $p$. We will call $\bar Q_\lambda$ the slow vector field, and its flow the slow dynamics. The time variable of the slow dynamics is the slow time $\tau=\epsilon t$ where $t$ is the time of $X_{\epsilon,\lambda}$. It is not difficult to see that
\begin{equation}\label{eq-SVF-example}
\bar Q_{\lambda}=-\frac{g(x,y,0,\lambda)\frac{\partial f}{\partial y}(x,y,0,\lambda)}{\frac{\partial f}{\partial x}(x,y,0,\lambda)}\frac{\partial}{\partial x}+g(x,y,0,\lambda)\frac{\partial}{\partial y}.
\end{equation}
Assumption 1 and L'Hospital's rule imply that the slow vector field $\bar Q_{\lambda_0}$ can be regularly extended through the slow-fast Hopf point $p_c$, with the $x$-component equal to 
$$-\frac{\frac{\partial g}{\partial x}(p_c,0,\lambda_0)\cdot \frac{\partial f}{\partial y}(p_c,0,\lambda_0)}{\frac{\partial^{2} f}{\partial x^{2}}(p_c,0,\lambda_0)}\ne 0$$
at $p_c$. Near $p_c$, the slow dynamics points from the normally attracting branch to the normally repelling branch of $\mathcal C_{\lambda_0}$.

Now, we define the notion of slow divergence integral \cite[Chapter 5]{DDR-book-SF}. Assume that $m_\lambda\subset \mathcal{C}_{\lambda}$ is a normally hyperbolic segment that does not contain singularities of the slow vector field $\bar Q_\lambda$, given by \eqref{eq-SVF-example}. The slow divergence integral associated with $m_\lambda$ is defined by
\begin{equation}
    \label{SDI-classicaldef1}
    I(m_\lambda)=\int_{\tau_1}^{ \tau_2}\divergenceOperator X_{0,\lambda}(z_\lambda(\tau))d \tau,
\end{equation}
where $$\divergenceOperator X_{0,\lambda}(x,y)=\frac{\partial f}{\partial x}(x,y,0,\lambda)$$ is the divergence of the fast subsystem $X_{0,\lambda}$, $z_\lambda:[\tau_1,\tau_2]\to\mathbb R^2$, $z_\lambda'(\tau)=\bar Q_\lambda(z_\lambda(\tau))$ and $z_\lambda(\tau_1)$ and $z_\lambda(\tau_2)$ are the endpoints of $m_\lambda$. Note that we parameterize $m_\lambda$ by the slow time $\tau$. This definition is independent of the choice of parameterization $z_\lambda$ of $m_\lambda$ (see \cite[Chapter 5]{DDR-book-SF}).

From the fact that $\bar Q_{\lambda_0}$ has a regular extension through $p_c$ it follows that we can compute  the slow divergence integral \eqref{SDI-classicaldef1} along any small segment $m_{\lambda_0}\subset \mathcal C_{\lambda_0}$, containing $p_c$ in its interior. More precisely, $m_{\lambda_0}$ is contained in a small neighborhood of $p_c$ and it consists of a normally attracting branch, a normally repelling branch and $p_c$ between them.
For $y\ne y_c$ and $y\sim y_c$, we denote by $\alpha(y)\in \mathcal C_{\lambda_0}$ (resp. $\omega(y)\in \mathcal C_{\lambda_0}$) the $\alpha$-limit point (resp. the $\omega$-limit point) of the fast orbit of $X_{0,\lambda_0}$ through $(x_c,y)$. If $\mathcal C_{\lambda_0}$ is concave up (resp. concave down), then we assume that $y>y_c$ (resp. $y< y_c$). Notice that the singularity $\alpha(y)$ is normally repelling and $\omega(y)$ is normally attracting. Finally, using \eqref{eq-SVF-example} and \eqref{SDI-classicaldef1} we define the slow divergence integral associated with $m_{\lambda_0}\subset \mathcal C_{\lambda_0}$ from $\omega(\tilde y)$ to $\alpha(\bar y)$, with $\tilde y,\bar y>y_c$ (resp. $\tilde y,\bar y< y_c$):
\begin{equation}
    \label{SDIFinal}
    I(\tilde y,\bar y)=-\int_{\omega(\tilde y)_x}^{\alpha(\bar y)_x}\frac{\left(\frac{\partial f}{\partial x}(x,y(x),0,\lambda_0)\right)^2}{g(x,y(x),0,\lambda_0)\frac{\partial f}{\partial y}(x,y(x),0,\lambda_0)}dx,
\end{equation}
where $\alpha(\bar y)_x$ is the $x$-coordinate of $\alpha(\bar y)$, $\omega(\tilde y)_x$ is the $x$-coordinate of $\omega(\tilde y)$ and $\mathcal C_{\lambda_0}$ is locally the graph of $y=y(x)$ ($y(x_c)=y_c$). We call $I(\tilde y,\bar y)=0$ the entry-exit function. 
 
We write $\tilde I( y)=I( y, y)$. Clearly, we have $\tilde I(y)\to 0$ as $ y\to y_c$.\\
\\
\textbf{Assumption 2} \textit{(Slow divergence integral).} Let $\mathcal C_{\lambda_0}$ be concave up (resp. concave down) near the slow-fast Hopf point $p_c$. We assume that $\tilde I( y)\ne 0$ for all $y\in J:=(y_c,y_c+\rho)$ (resp. $y\in J=:(y_c-\rho,y_c)$), with $\rho>0$ small enough. \\
\\
Let $y_0\in J$ be fixed. From Assumptions 1 and 2 and \eqref{SDIFinal} it follows that the entry-exit function $I(\tilde y,\bar y)=0$ produces a sequence $(y_k)_{k\ge 0}$ in $J$ that converges monotonically to $y_c$. More precisely, 
\begin{itemize}
    \item ($\tilde I( y)< 0$ for all $y\in J$) In this case, the sequence $(y_k)_{k\ge 0}$ is generated by $I(y_{k+1},y_{k})=0$, $k\ge 0$, and tends monotonically to $y_c$ as $k\to\infty$. 
    
    \item ($\tilde I( y)> 0$ for all $y\in J$) In this case, we use $I(y_k,y_{k+1})=0$, $k\ge 0$, and $(y_k)_{k\ge 0}$ tends monotonically to $y_c$ as $k\to\infty$. 
\end{itemize}
We call $(y_k)_{k\ge 0}$ the entry-exit sequence starting at $y_0\in J$.
When $\mathcal C_{\lambda_0}$ is concave up (resp. concave down), the sequence $(y_k)_{k\ge 0}$ is decreasing (resp. increasing). See, e.g., Figure \ref{fig-Hopf}.

\subsection{Fractal classification of the slow-fast Hopf point}\label{subsection-Hopf}
The cyclicity of $p_c$ in $X_{\epsilon,\lambda}$ is bounded by $M\in\mathbb{N}_0$ if there exist $\epsilon_0>0$, a neighborhood $V$ of $\lambda_0$ and a neighborhood $W$ of $p_c$ such that $X_{\epsilon,\lambda}$ has at most $M$ limit cycles in $W$ for each $(\epsilon,\lambda)\in [0,\epsilon_0]\times V$. The smallest $M$ with this property is called the cyclicity of $p_c$ in the slow-fast family $X_{\epsilon,\lambda}$. We denote it by $\cycl(X_{\epsilon,\lambda},p_c)$.

Theorem \ref{theorem-Main1} below follows from \cite[Theorem 3.2]{BoxNovo} and \cite[Theorem 3.4]{BoxNovo} (see also \cite{BoxVlatko}).
\begin{theorem}
\label{theorem-Main1}
Consider a smooth slow-fast system $X_{\epsilon,\lambda}$ of the form in \eqref{eq-def-slowfast-1} that satisfies Assumptions 1 and 2. Let $(y_k)_{k\ge 0}$ be an entry-exit sequence defined in Section \ref{subsection-SDI}. Then the following statements are true. 
\begin{enumerate}
    \item The Minkowski dimension $\dim_B (y_k)_{k\ge 0}$ exists and \begin{equation}\label{equations-bijection}
  \dim_B\mathcal (y_k)_{k\ge 0}\in \left\{\frac{2j+1}{2j+3} \ | \ j=0,1,2,\dots \right\}\cup \{1\}.  \nonumber
\end{equation} 
If $\dim_B (y_k)_{k\ge 0}\ne 1$, then $ (y_k)_{k\ge 0}$ is Minkowski nondegenerate. Furthermore, the Minkowski dimension and the Minkowski nondegeneracy of $ (y_k)_{k\ge 0}$  do not depend on the choice of the first element $y_0\in J$ of $(y_k)_{k\ge 0}$.
\item If $X_{\epsilon,\lambda}$ is analytic, then $\cycl(X_{\epsilon,\lambda},p_c)$ is finite. Moreover, if $\dim_B (y_k)_{k\ge 0}\ne 1$, then $$\cycl(X_{\epsilon,\lambda},p_c)\le \frac{\dim_B (y_k)_{k\ge 0}+1}{2\left(1-\dim_B (y_k)_{k\ge 0}\right)}.$$
\end{enumerate}
\end{theorem}
Let us briefly explain the idea behind the proof of Theorem \ref{theorem-Main1}. If $\lambda=\lambda_0$ is fixed, then we can bring the system $X_{\epsilon,\lambda_0}$, locally near $p_c$, into  the following normal
form for $C^\infty$-equivalence:
\begin{equation}
\label{normal form 2}
     \begin{vf}
        \dot{x} &=& y-x^{2}  \\
        \dot{y} &=&\epsilon\left( g_1(x,\epsilon)+\left(y-x^2\right)g_2(x,y,\epsilon)\right),
    \end{vf}
\end{equation}
where $g_1$ and $g_2$ are smooth functions and $g_1(x,0)=-x+ O(x^2)$. ($C^\infty$-equivalence means $C^\infty$-coordinate change and a multiplication by
a $C^\infty$-positive function.) The slow-fast Hopf point of \eqref{normal form 2} is located at the origin $(x,y)=(0,0)$. We can give a complete fractal classification of entry-exit sequences associated with the normal form \eqref{normal form 2} (see \cite[Theorem 3.1]{BoxNovo} and \cite{EZZ}). Then Theorem \ref{theorem-Main1}.1 follows from the invariance of the slow divergence integral under $C^\infty$-equivalences \cite[Chapter 5]{DDR-book-SF} and the invariance of the Minkowski dimension under bi-Lipschitz maps (see Section \ref{subsection-Mink}). Theorem \ref{theorem-Main1}.2 is a consequence of \cite[Lemma 3.3]{BoxNovo}, \cite{DRbirth} and \cite{HuzakNormal}. For a detailed proof we refer to \cite{BoxNovo}.

\begin{remark}\label{remark-important}
Theorem \ref{theorem-Main1}.1 implies that $\dim_B (y_k)_{k\ge 0}$ can take only the following discrete set of values: $1/3$, $3/5$, $5/7,\dots,1$. In the analytic case, if the Minkowski dimension is $1/3$, then the slow-fast Hopf point $p_c$ can produce at most $1$ limit cycle; if the Minkowski dimension is $3/5$, then $p_c$ produces at most $2$ limit cycles, etc. See Theorem \ref{theorem-Main1}.2. In \cite{BoxNovo}, we call the upper bound for $\cycl(X_{\epsilon,\lambda},p_c)$ given in Theorem \ref{theorem-Main1}.2 the fractal codimension of $p_c$. 
\end{remark}
We can numerically calculate the Minkowski dimension $\dim_B (y_k)_{k\ge 0}$ using explicit fractal formulas of Cahen-type (see \cite[Proposition 1]{BoxVesna}). Then we find the upper bound for $\cycl(X_{\epsilon,\lambda},p_c)$ (Theorem \ref{theorem-Main1}.2). We refer the reader to \cite[Section 6]{BoxNovo} for some interesting numerical examples (a two-stroke oscillator, classical Li\'{e}nard equations, $\dots$).

\subsection{Fractal classification of canard cycles}\label{subsection-canard}
In this section, we consider a smooth slow-fast system $X_{\epsilon,\lambda}$ of the form in \eqref{eq-def-slowfast-1} that satisfies Assumption 1. 

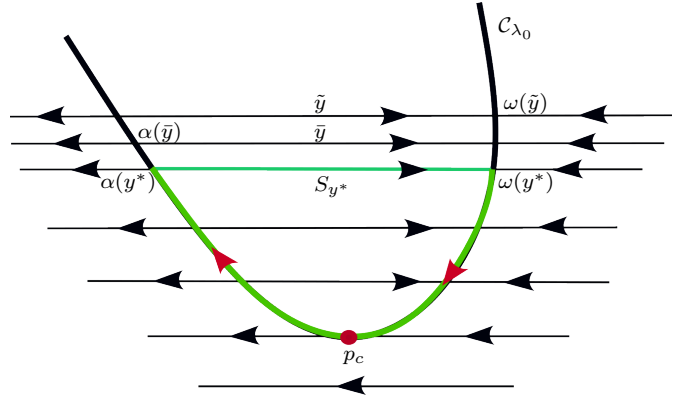
\begin{figure}[htb]
	\begin{center}
		\includegraphics[width=8.7cm,height=5.2cm]{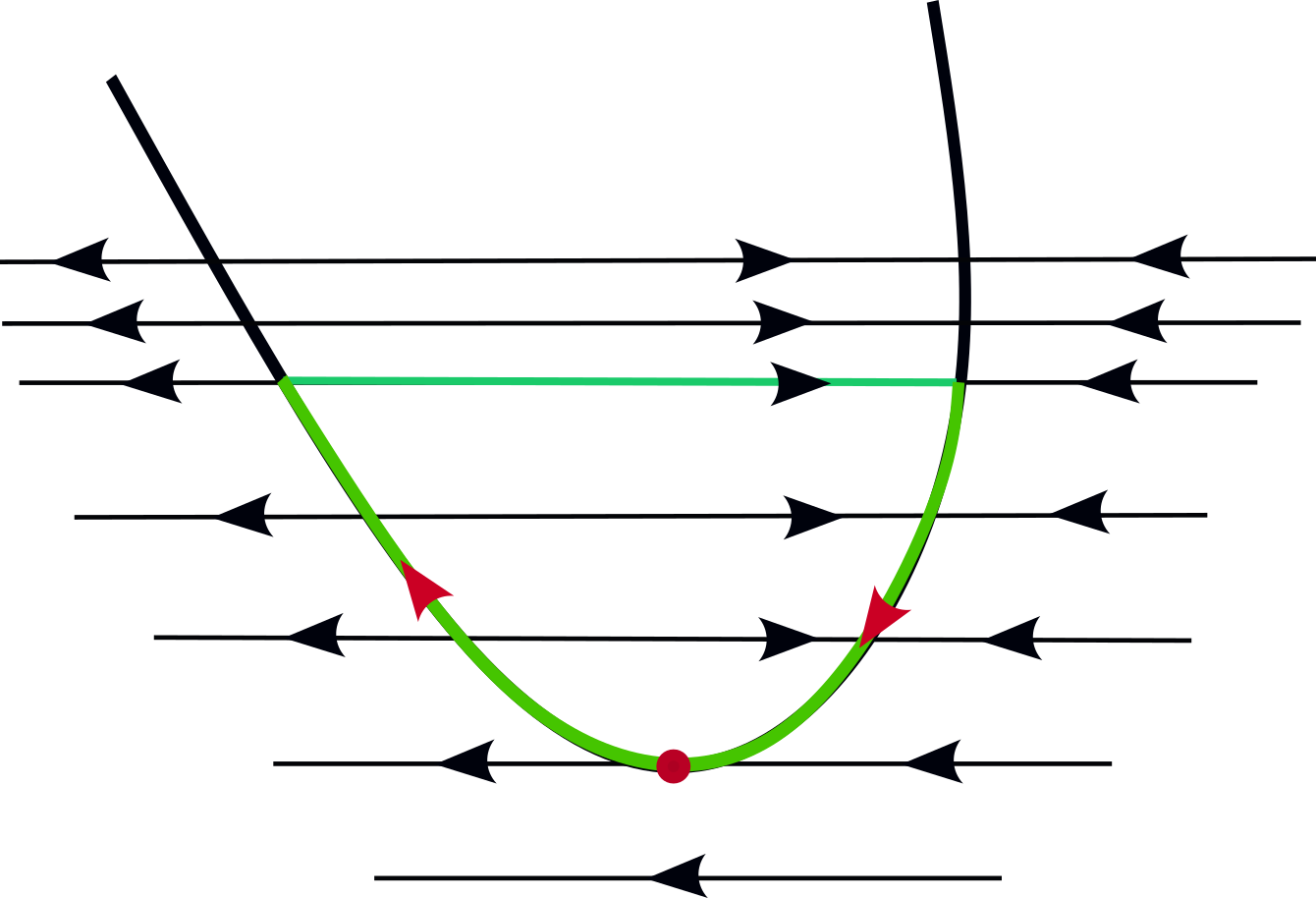}
		{\footnotesize \put(-123,13){$p_c$}
        \put(-65,137){$\mathcal C_{\lambda_0}$}
        \put(-65,78){$\omega(y^*)$}
        \put(-63,108){$\omega(\tilde y)$}
        \put(-200,97){$\alpha(\bar y)$}
        \put(-214,78){$\alpha(y^*)$}
        \put(-134,77){$S_{y^*}$}
        \put(-134,108){$\tilde y$}
        \put(-134,97){$\bar y$}
}
         \end{center}
	\caption{Dynamics of $X_{0,\lambda_0}$ and the canard cycle $S_{y^*}$. $\mathcal C_{\lambda_0}$ is concave up near $p_c$.}
	\label{fig-Canard}
\end{figure}

Let $\mathcal C_{\lambda_0}$ be concave up (resp. concave down) near the slow-fast Hopf point $p_c=(x_c,y_c)$. We take $y^*>y_c$ (resp. $y^*<y_c$) and define the canard cycle $S_{y^*}$ at the level $(\epsilon,\lambda)= (0, \lambda_0)$: It consists of a fast orbit of $X_{0,\lambda_0}$ located inside $y=y^*$ and the
portion $m_{\lambda_0}$ of the curve of singularities $\mathcal C_{\lambda_0}$ between the $\alpha$-limit point $\alpha(y^*)\in \mathcal C_{\lambda_0}$ and the $\omega$-limit point $\omega(y^*)\in \mathcal C_{\lambda_0}$ of that orbit. We assume that the set $m_{\lambda_0}\setminus \{p_c\}$ is normally hyperbolic. Then the slow vector field $\bar Q_{\lambda_0}$ is well defined along the segment $m_{\lambda_0}$ (see Section \ref{subsection-SDI}). If $\bar Q_{\lambda_0}$ is regular along $m_{\lambda_0}$ (pointing from the attracting branch to the repelling branch of $\mathcal C_{\lambda_0}$), then the slow divergence integral in \eqref{SDI-classicaldef1} is well defined along the segment $m_{\lambda_0}$. More precisely, we can again consider the slow divergence integral $I(\tilde y,\bar y)$ associated with the segment in $\mathcal C_{\lambda_0}$ from $\omega(\tilde y)$ to $\alpha(\bar y)$, with $\tilde y,\bar y\sim y^*$. We refer to Figure \ref{fig-Canard}.

We write $\tilde I( y)=I( y, y)$, for $y\sim y^*$.  
We assume that $S_{y^*}$ is balanced (that is, $\tilde I( y^*)=0$). 
Instead of Assumption 2 (see Section \ref{subsection-SDI}) we need 
\\
\\
\textbf{Assumption 2'} \textit{(Slow divergence integral).} Let $\mathcal C_{\lambda_0}$ be concave up (resp. concave down) near the slow-fast Hopf point $p_c$. We assume that $\tilde I( y)\ne 0$ for all $y\in J:=(y^*,y^*+\rho)$ (resp. $y\in J:=(y^*-\rho,y^*)$), with $\rho>0$ small enough. \\
\\
Now we can define, in a similar fashion as before, entry-exit sequences associated with the canard cycle $S_{y^*}$. 
Let $y_0\in J$ be fixed. We consider the following two cases depending on whether $\tilde I$ is negative or positive in $J$: 
\begin{itemize}
    \item ($\tilde I( y)< 0$ for all $y\in J$) The entry-exit sequence $(y_k)_{k\ge 0}$ is generated by $I(y_{k+1},y_{k})=0$, $k\ge 0$, and tends monotonically to $y^*$ as $k\to\infty$. 
    
    \item ($\tilde I( y)> 0$ for all $y\in J$) In this case, we use $I(y_k,y_{k+1})=0$, $k\ge 0$, and $(y_k)_{k\ge 0}$ tends monotonically to $y^*$ as $k\to\infty$. 
\end{itemize}
If $\mathcal C_{\lambda_0}$ is concave up (resp. concave down) near $p_c$, the sequence $(y_k)_{k\ge 0}$ is decreasing (resp. increasing).

Now, we have (see \cite{BoxRenato} or \cite[Theorem 3]{BoxVlatko})
\begin{theorem}\label{theorem-Main2}
  Consider a smooth slow-fast system $X_{\epsilon,\lambda}$ of the form in \eqref{eq-def-slowfast-1} that satisfies Assumptions 1 and 2'. Let $(y_k)_{k\ge 0}$ be an entry-exit sequence defined above. Then the following statements are true. 
\begin{enumerate}
    \item The Minkowski dimension $\dim_B (y_k)_{k\ge 0}$ exists and \begin{equation}
  \dim_B\mathcal (y_k)_{k\ge 0}\in \left\{\frac{j}{j+1} \ | \ j=0,1,2,\dots \right\}\cup \{1\}.  \nonumber
\end{equation} 
If $\dim_B (y_k)_{k\ge 0}\ne 0,1$, then $ (y_k)_{k\ge 0}$ is Minkowski nondegenerate. Furthermore, the Minkowski dimension and the Minkowski nondegeneracy of $ (y_k)_{k\ge 0}$ do not depend on the choice of the first element $y_0$ of $(y_k)_{k\ge 0}$.
\item If $\dim_B (y_k)_{k\ge 0}\ne 1$, then $$\cycl(X_{\epsilon,\lambda},S_{y^*})\le \frac{2-\dim_B (y_k)_{k\ge 0}}{1-\dim_B (y_k)_{k\ge 0}}.$$
\end{enumerate}
\end{theorem}

We denote by $\cycl(X_{\epsilon,\lambda},S_{y^*})$ the cyclicity of $S_{y^*}$ in the family $X_{\epsilon,\lambda}$ (see Theorem \ref{theorem-Main2}.2). The definition of $\cycl(X_{\epsilon,\lambda},S_{y^*})$ is similar to the definition of
$\cycl(X_{\epsilon,\lambda},p_c)$ given in Section \ref{subsection-Hopf} ($W$ is now a Hausdorff neighborhood of $S_{y^*}$).

Theorem \ref{theorem-Main2} establishes a bijective correspondence between the Minkowski dimension $\dim_B (y_k)_{k\ge 0}$ ($0$, $1/2$, $2/3,\dots$) and the upper bound for the cyclicity of $S_{y^*}$ ($2$, $3$, $4,\dots$).

	\label{fig-chirp1}

\section{Conclusion and discussions}
In this paper, the Minkowski dimension is used to measure the local density of spiral trajectories or countable sets accumulating at a fixed point. The bigger the Minkowski dimension, the more limit cycles can be born after perturbation. In the regular case, the fractal analysis of orbits generated by Poincar\'{e} map is a useful tool to find out how many limit cycles can be produced by monodromic limit periodic sets. In the slow-fast case, the Poincar\'{e} map is not defined in the limit $\epsilon\to 0$, and it is natural to generate countable sets accumulating at a fixed point using the concept of a slow divergence integral and entry-exit function. The density of such sets plays a crucial role when we deal with the birth of canards. 

A fractal analysis of focus-like fixed points of piecewise smooth systems is an open direction for further research. Piecewise smooth systems are often used to model electrical circuits that have switches, problems with friction, control systems, etc. We refer to \cite{PWS}. We believe that the Minkowski dimension of spiral trajectories near such focus-like fixed points is related to the number of crossing limit cycles \cite{Rom} after perturbation. 

Another open direction for further research is a fractal analysis of slow-fast systems in $\mathbb R^3$ with two fast variables and one slow variable. The entry-exit functions \cite{Popovic2025}, together with the presence of a return mechanism, can produce closed trajectories. Again, we believe that the Minkowski dimension is related to the creation of such closed trajectories.

Fractal analysis of spiral trajectories in homoclinic bifurcations in three-dimensional systems is also an interesting problem (see Shil’nikov theorems in the saddle-focus case \cite[Section 6.3]{Kuz}). Shil’nikov’s results can be applied to the modeling of nerve impulses.

\section*{Declarations}
 
\textbf{Ethical Approval} \ 
Not applicable.
 \\
\\
 \textbf{Competing interests} \  
The authors declare that they have no conflict of interest.\\
 \\
\textbf{Funding} \
The research of R. Huzak, G. Radunovi\'{c} and V.\ \v Zupanovi\'c was supported by: Croatian Science Foundation (HRZZ) grant IP-2022-10-9820. Additionally, the research of G. Radunovi\'{c} and V.\ \v Zupanovi\'c was partially
supported by the Horizon grant 101183111-DSYREKI-HORIZON-MSCA-2023-SE-01.\\
 \\
\textbf{Availability of data and materials}  \
Not applicable.

\bibliographystyle{plain}
\bibliography{bibtex}
\end{document}